\documentclass[nodea]{birkmult}
\usepackage{amsfonts,amsmath,amsthm,amscd,amssymb,latexsym}
\usepackage{color}

\newtheorem{thm}{Theorem}[section]
 \newtheorem{cor}[thm]{Corollary}
 
 \newtheorem{prop}[thm]{Proposition}
 \theoremstyle{definition}
 \newtheorem{defn}[thm]{Definition}
 \theoremstyle{remark}
 \newtheorem{rem}[thm]{Remark}
 
 \numberwithin{equation}{section}

\def\dom{{\mathfrak D}}
\newcommand\diag{\operatorname{diag}}
\newcommand{\ctg}{{\mathrm{ctg}}}
\newcommand{\spa}{{\mathrm{span}}}

\DeclareMathOperator*{\rank}{rank}
\DeclareMathOperator*{\ran}{ran} \DeclareMathOperator*{\sgn}{sgn}

\def\gl{{\lambda}}
\def\gG{{\Gamma}}

\def\gotH{{\mathfrak H}}

\def\cH{{\mathcal H}}

\def\dC{{\mathbb C}}
\def\dR{{\mathbb R}}
\def\dN{{\mathbb N}}
\def\dZ{{\mathbb Z}}

\def\kH{{\mathcal H}}

\begin{document}
\title[1-D
Schr\"odinger Operators  with  Point Interactions]
 {On the Negative Spectrum  of  One-Dimensional
Schr\"odinger Operators  with Point Interactions}

\author{N.~Goloschapova}

\address{%
Institute of Applied Mathematics and Mechanics, NAS of Ukraine,\\
R. Luxemburg str. 74\\
83114 Donetsk,Ukraine}

\email{ng85@bk.ru}
\author{L.~Oridoroga}
\address{%
Donetsk National University\\
Universitetskaja str. 24\\
83055 Donetsk,
Ukraine} \email{oridoroga@skif.net}



\date{January 27, 2009}

\begin{abstract}
We investigate negative spectra of 1--D Schr\"odinger operators
with $\delta$- and $\delta'$-interactions on a discrete set in the
 framework  of a new approach. Namely, using  technique of boundary triplets
 and the
corresponding Weyl functions, we complete and generalize the
results of S. Albeverio and L. Nizhnik \cite{AlbNiz03, AlbNiz03a}.
For instance,  we propose  the algorithm for determining the
number of negative  squares  of  the operator  with
$\delta$-interactions. We also  show that the number of negative
squares of the operator with $\delta'$-interactions equals the
number of negative strengths.
\end{abstract}

\maketitle
\section{Introduction}
Consider formal differential expressions
\begin{equation}\label{eq0}
              \ell_{X,\alpha}=-\frac{\mathrm{d}^2}{\mathrm{d}x^2}+\sum^{\infty}_{k=1}\alpha_k\delta_k(x),\quad 
              \ell_{X,\beta}=-\frac{\mathrm{d}^2}{\mathrm{d}x^2}+\sum^{\infty}_{k=-\infty}\beta_k\langle\cdot,\delta_k'\rangle\delta_k'(x),
       \end{equation}
where $\delta_k(x):=\delta(x-x_k)$ and $\delta(x)$ is a Dirack
delta-function, $\alpha_k,\beta_k\in\dR$, $k\in I$, and $I$ equals either $\dN$ or $\dZ$. We assume that $X=\{x_k\}_{k\in I}\subset\mathbb{R}$ is an
increasing sequence such that $d_k:= x_{k+1}-x_k>0$, $k\in I$, and
\begin{equation}\label{d_*^*}
 d_*:=\inf_{k\in I}d_k>0,\,\qquad d^*:=\sup_{k\in I}d_k<\infty.
\end{equation}
One defines the corresponding operators with $\delta$-
 and $\delta'$-interactions  in
$L^2(\dR)$ by the differential expression
$-\frac{\mathrm{d}^2}{\mathrm{d}x^2}$ on the domains, respectively,
\begin{gather}\label{delta_op}
 \dom(L_{X,\alpha})=\bigl\{f\in
W_2^2(\mathbb{R}\setminus X):\begin{array}{c}    f(x_k+)=f(x_k-), \\
                                                             f'(x_k+)-f'(x_k-)=\alpha_k f(x_k) \\
                                                          \end{array},\,
                                                          x_k\in X\bigr\},\\
\dom(L_{X,\beta})=\bigl\{f\in
W_2^2(\mathbb{R}\setminus X):\begin{array}{c}    f'(x_k+)=f'(x_k-), \\
                                                             f(x_k+)-f(x_k-)=\beta_kf'(x_k) \\
                                                          \end{array},\,
                                                          x_k\in X\bigr\}.\label{delta'_op}
\end{gather}
Note that the operators $L_{X,\alpha}$ and
$L_{X,\beta}$ are self-adjoint (\cite{AlbGes88}, see also \cite{bsw,gk}).


Schr\"odinger operators with point interactions have been studied
extensively in the last decades (numerous results and a
comprehensive list of references may be found in \cite{AlbGes88,
AlbKur00}, see also Appendix K by P. Exner in \cite{AlbGes88}). In
the recent publications \cite{AlbNiz03, AlbNiz03a}, S. Albeverio
and L. P. Nizhnik investigated the numbers
$\kappa_-(L_{X,\alpha})$ and $\kappa_-(L_{X,\beta})$ of negative
eigenvalues  of the operators $L_{X,\alpha}$ and $L_{X,\beta}$
 in the case $|X|=n<\infty$. They described $\kappa_-(L_{X,\alpha})$ in terms of a certain continued
  fractions (cf. \cite[Theorem 3]{AlbNiz03}) and also proposed an elegant algorithm
for determining $\kappa_-(L_{X,\alpha})$. In particular, they
formulated necessary and sufficient
 conditions  in terms of the distances $d_k$ and the strengths $\alpha_k$ for the equalities $\kappa_-(L_{X,\alpha})=n$
 and $\kappa_-(L_{X,\alpha})=0$ to hold (cf. \cite[Theorem 5]{AlbNiz03} and \cite[Theorem 4]{AlbNiz03}, respectively).
Regarding the operators with $\delta'$-interactions, it is shown
in \cite[Theorem 6]{AlbNiz03a} that the number of negative
eigenvalues of  $L_{X,\beta}$ equals $n$ if and only if all
intensities are negative, i.e., $\kappa_-(\{\beta_k\}_{k=1}^n)=n$.

  In this  paper, we present a new approach to investigate negative spectra of the operators with $\delta$-
  and $\delta'$-interactions
    on the discrete  set $X$ satisfying \eqref{d_*^*}. Namely, we consider the operators $L_{X,\alpha}$ and $L_{X,\beta}$ as self--adjoint extensions of the symmetric operator
\begin{equation}\label{eq26}
L_{\min}= -\frac{\mathrm{d}^2}{\mathrm{d}x^2},\qquad
\dom(L_{\min})=\text{\r{W}}_{2}^{2}(\mathbb{R}\setminus X),\qquad
X=\{x_k\}_{k\in I}
\end{equation}
and apply the technique  of boundary  triplets and the
corresponding Weyl functions (see  \cite{Gor84, DerMal91} and also
Section \ref{prelim}). We construct a boundary triplet for
$L_{min}^*$ and establish a connection between  Hamiltonians
$L_{X,\alpha}$  and $L_{X,\beta}$ and a certain classes  of Jacobi
matrices. Using this  connection, we describe
 $\kappa_-(L_{X,\alpha})$ and $\kappa_-(L_{X,\beta})$ by  means of entries of these matrices (Theorem \ref{th12}).
  The latter enables us to complete and
substantially generalize previous results from \cite{AlbNiz03,
AlbNiz03a} mentioned above.
 Namely, for a $\delta$-type interactions, we construct an
algorithm for determining $\kappa_-(L_{X,\alpha})$ (Theorem
\ref{th2}). In the case $|X|=n$, our algorithm differs from the
one proposed  by S. Albeverio and L.P. Nizhnik, but it  is  close
to that (see Remark \ref{AN appr}).
One of our  main  results  is the following  equality
$\kappa_-(L_{X,\beta})=\kappa_-(\beta)$ (Theorem 4.1). It means
that the number of negative  squares  of  $L_{X,\beta}$ equals the
number of negative intensities.
In the particular case $\kappa_-(\beta) =|X|= n < \infty$, this
results coincides with \cite[Theorem 6]{AlbNiz03a}. It is
interesting to mention that for   the operator with
$\delta$-interactions   such equality does  not hold (cf.
\cite{AlbNiz03,Ogu08}). We obtain sufficient  condition for the
inequality $\kappa_-(L_{X,\alpha})\geq m$ (as well as for the
equality)  with any $m$ (Theorem \ref{suffcond}). It differs from
the one recently obtained by Ogurisu in \cite{Ogu08} and implies
sufficient condition for $\kappa_-(L_{X,\alpha})=n$ proposed  by
Albeverio and Nizhnik \cite[Criterion 3]{AlbNiz03a} in the case
$\kappa_-(\alpha)=|X|=n$. In particular, the operator
$L_{X,\alpha}$ with arbitrary number of negative intensities might
be non-negative (see \cite[Theorem 4]{AlbNiz03} and also Corollary
\ref{remunus}).

The results of the paper  were partially  announced (without
proofs) in \cite{Gol08}.

\textbf{Notation.} 
Let $X$ be a discrete subset of $\mathbb{R}$; $|X|$ stands for the
cardinal number  of the set $X$. By $W^2_2(\mathbb{R}\setminus X)$
and
$\text{\r{W}}_{2}^{2}(\mathbb{R}\setminus X)$ we denote the Sobolev spaces
\begin{gather*}
W^2_2(\mathbb{R}\setminus X):=\{f\in L^2(\mathbb{R}):\ f, f'\in
AC_{loc}(\mathbb{R}\setminus X),\  f''\in
  L^2(\mathbb{R})\},\\
\text{\r{W}}_{2}^{2}(\mathbb{R}\setminus X):=\{f\in
W^2_2(\mathbb{R}): f(x_k)=f'(x_k)=0\,\  \mbox{for all }  x_k\in
X\}.
\end{gather*}

\section{Preliminaries}\label{prelim}

\textbf{Boundary triplets and closed extensions.} In this
subsection, we recall basic notions of the theory of boundary
triplets (we refer the reader to \cite{DerMal91, Gor84} for a
detailed exposition).

Let $A$ be a closed densely defined symmetric operator in the
Hilbert space $\gotH$ with equal deficiency indices
$n_\pm(A)=\dim\ker(A^*\pm i)\leq\infty$.
\begin{defn}[\cite{Gor84}]\label{bound}%
A triplet $\Pi=\{\kH,\gG_0,\gG_1\}$ is called a {\rm boundary
triplet} for the adjoint operator $A^*$ of $A$ if $\kH$ is an auxiliary
Hilbert space and
$\Gamma_0,\Gamma_1:\  \dom(A^*)\rightarrow\kH$ are linear mappings such that\\
 $(i)$ the  second  Green identity,
\begin{equation*}
(A^*f,g)_\gotH - (f,A^*g)_\gotH = (\gG_1f,\gG_0g)_\kH - (\gG_0f,\gG_1g)_\kH,
\end{equation*}
holds for all $f,g\in\dom(A^*)$, and\\
$(ii)$ the mapping $\gG:=(\Gamma_0,\Gamma_1)^\top: \dom(A^*)
\rightarrow \kH \oplus\kH$ is surjective.
\end{defn}
Since $n_+(A)=n_-(A)$, a boundary triplet
$\Pi=\{\kH,\gG_0,\gG_1\}$ for $A^*$ exists and is not unique
\cite{Gor84}. Moreover, $\dim\kH=n_\pm(A) $ and
$A=A^*\upharpoonright\ker(\Gamma_0)\cap\ker(\Gamma_1)$.
Any  proper extension $\widetilde{A}$ of $A$
admits the following representation  (see \cite{dhms})
\begin{equation}\label{eq80}
\dom(\widetilde{A})=\dom(A_{C,D}):=\dom(A^*)\upharpoonright\ker(D\Gamma_1-C\Gamma_0),\quad\text{where
}\quad C,\ D\in[\kH].
\end{equation}
Note that representation \eqref{eq80} is not unique.

In what follows, we  will also denote
$A_0=A^*\upharpoonright\ker(\Gamma_0)$. Note that $A_0^*=A_0$.

\begin{defn}[\cite{DerMal91}]\label{Weylfunc}
Let $\Pi=\{\cH,\gG_0,\gG_1\}$ be a boundary triplet  for $A^*.$
The operator valued function  $M(\cdot) :\rho(A_0)\rightarrow
[\cH]$ defined by
\begin{equation}\label{2.3A}
\Gamma_1 f_\gl=M(\gl)\Gamma_0f_\gl, \qquad \gl\in\rho(A_0),\quad
f_\gl\in \mathcal{N}_\gl=\ker(A^*-\gl),
      \end{equation}
is called  {\em the Weyl function}  corresponding to the boundary
triplet $\Pi.$
\end{defn}
%

Before formulate next result we need the
following definition.
\begin{defn} [\cite{Kat}]
 Let   $T=T^*\in \mathcal{C}(\gotH)$ 
 and let $E_T(\lambda)=E_T(\lambda-0)$ be the spectral function of $T$.
     Dimension of the subspace $E_T(-\infty,0)\gotH$ is called a number of
     negative squares of $T$ and is denoted by  $\kappa_{-}(T)$.
     \end{defn}

   The Weyl function $M(\cdot)$ enables us to describe the number of negative
squares of self-adjoint extensions of $A$.
  \begin{thm}[\cite{DerMal91}]\label{th8}
 Let $A$ be nonnegative, $A\geq 0$, and let $A_{C,D}$ be its self-adjoint
extension. Assume that $\Pi=\{\mathcal{H}, \Gamma_0, \Gamma_1\}$ is a boundary
triplet for $A^*$ such that $A_0=A_F$, where $A_F$ is the Friedrichs extension of $A$.
 If the strong resolvent limit
$M(0):=s-R-\lim\limits_{x\uparrow 0}M(x)$ (see \cite[Capter 8]{Kat})
 exists and $M(0)\in [H]$, then
\begin{equation}\label{eq40'}
\kappa_{-}(A_{C,D})= \kappa_{-}(CD^*-D M(0)D^*).
\end{equation}
 \end{thm}
\quad \\
\textbf{The Sylvester criterion.}
Description  of $\kappa_-(L_{X,\alpha})$ is substantially based on
 the  following  fact  (see, for instance, \cite[Lemma
4]{Mal}).
 \begin{prop}\label{th2'} Let the operator
$T=T^*\in\mathcal{C}(\mathcal{H})$ admit the block-matrix
representation
$
T=\left(%
\begin{array}{cc}
  T_{11} & T_{12} \\
  T_{21} & T_{22} \\
\end{array}%
\right)
$ 
with  respect to the decomposition
$\mathcal{H}=\mathcal{H}_1\oplus \mathcal{H}_2$, where
$T_{11}\in[{\mathcal{H}}_1]$,
$T_{12}=T_{12}^*\in[\mathcal{H}_2,\mathcal{H}_1]$ and
$T_{22}\in\mathcal{C}(\mathcal{H}_2)$. If $0 \in \rho(T_{11})$,
then
\begin{equation}
\kappa_{-}(T)=
\kappa_{-}( T_{11}) + 
\kappa_{-}(T_{22}-T_{21}T_{11}^{-1}T_{12}).
\end{equation}
\end{prop}
\section{Operators with $\delta$ -- type  interactions}
\subsection{The case of infinite number of $\delta$ -- type
interactions}\label{infinite}

       Consider the following Jacobi matrix in  $l^2(\mathbb{N})$,
     \begin{equation}\label{eq5}
S=
\left(%
\begin{array}{cccc}
  \alpha_1+d^{-1}_1 & -d^{-1}_1 & 0 &\ldots  \\
  -d^{-1}_1 & \alpha_2+d^{-1}_1+d^{-1}_2 & -d^{-1}_2 & \ldots  \\
  0 & -d^{-1}_2 & \alpha_3+d^{-1}_2+d^{-1}_3   &\ldots  \\
  \ldots & \ldots &\ldots  &\ldots  \\
\end{array}%
\right).
\end{equation}
Notice that $S=S^*$ since $d_*=\inf_{k\in\mathbb{N}} d_k>0$ (cf.
\cite[Theorem VII.1.5]{Ber}).

      The main result of this  Section is the following description of  $\kappa_{-}(L_{X,\alpha})$.
\begin{thm}\label{th12}
Let the set $X=\{x_k\}_{k=1}^\infty$ satisfy \eqref{d_*^*}. Let also the operator  $L_{X,\alpha}$  and  the matrix $S$  be defined by \eqref{delta_op} and \eqref{eq5}, respectively.
Then $\kappa_{-}(L_{X,\alpha})=\kappa_{-}(S).$
\end{thm}
\begin{proof}
Consider the minimal operator \eqref{eq26}.
 Note that $n_{\pm}(L_{\min})=\infty$.
  Since $X$ satisfies \eqref{d_*^*}, the totality  $\Pi=\{\kH,\Gamma_0,\Gamma_1\}$, where
\begin{gather}\label{eq27}
\mathcal{H}=\oplus^{\infty}_{k=0}\mathcal{H}_k ,\qquad
\Gamma_0=\oplus^{\infty}_{k=0}\Gamma_0^k,\qquad
\Gamma_1=\oplus^{\infty}_{k=0}\Gamma_1^k,\quad \\
\label{eq27'}
\mathcal{H}_0=\mathbb{C},\qquad\Gamma_0^0f=-f(x_1-),\qquad\Gamma_1^0f=f'(x_1-),\quad \text{and}\\
\label{eq27''}
\mathcal{H}_k=\mathbb{C}^2,\qquad\Gamma_0^kf=\binom{f(x_k+)}{-f(x_{k+1}-)},
\qquad\Gamma_1^kf=\binom{f'(x_k+)}{f'(x_{k+1}-)},\quad k\in\dN,
 \end{gather}
  forms a boundary triplet for  $L^*_{\min}$  \cite[Lemma 1]{Koch}.
  The corresponding Weyl function is
\begin{gather}
M(\lambda)=\oplus_{k=0}^\infty M_k(\lambda),\quad
M_0(\lambda)=i\sqrt{\lambda},\\
 M_k(\lambda)=
\left(%
\begin{array}{cc}
  -\sqrt{\lambda}\ \ctg(\sqrt{\lambda}d_k) & -\sqrt{\lambda}/\sin(\sqrt{\lambda}d_k) \\
 -\sqrt{\lambda}/\sin(\sqrt{\lambda}d_k) & -\sqrt{\lambda}\ \ctg(\sqrt{\lambda}d_k) \\
\end{array}%
\right),\quad k\in\dN.\label{weyl}
\end{gather}
Using  \eqref{eq27}--\eqref{eq27''}, we obtain the representation
\eqref{eq80} for $\dom (L_{X,\alpha})$, where
\begin{equation*}
C=
 \left(%
   \begin{array}{ccccc}
   0 & \alpha_1 & 0 &0 &\ldots   \\
   1 & 1 & 0 & 0 &  \ldots   \\
  0 & 0 & 0 & \alpha_2 &  \ldots \\
  0 & 0 & 1 & 1 &  \ldots   \\
\ldots&\ldots&\ldots& \ldots&\ldots   \\
  \end{array}%
\right),\qquad
D=\left(%
   \begin{array}{cccccc}
  -1 & 1 & 0 &0 & \ldots     \\
   0&0&0  &0  & \ldots    \\
0 & 0 & -1 & 1 &    \ldots \\
0&0&0&0&\ldots \\
 \ldots &\ldots &  \ldots &\ldots& \ldots  \\
  \end{array}%
\right).
\end{equation*}
Without loss of generality, it can be assumed that  $C$ is
bounded. If there exists $\{\alpha_{k_j}\}_{j=1}^\infty$ such that
$\lim_{j\to\infty}|\alpha_{k_j}|=\infty$, then we put
\begin{gather*}
\widetilde{C}=KC,\quad \widetilde{D}=KD,\quad
\text{and}\quad \dom(L_{X,\alpha})=
\dom(L_{\min}^*)\upharpoonright\ker(\widetilde{D}\Gamma_1-\widetilde{C}\Gamma_0),\\
 \text{ where }
K=\diag(1,...,1,\alpha_{k_1}^{-1},1,...,1,\alpha_{k_2}^{-1},1,..).\qquad\qquad\qquad\qquad\qquad
\end{gather*}
 After straightforward calculations we get the matrix
 $T:=CD^*-DM(0)D^*$,
\begin{equation*}
T=\left(%
\begin{array}{cccccc}
  \alpha_1+d^{-1}_1 & 0 & -d^{-1}_1 & 0 & 0&\ldots  \\
  0 & 0 & 0 & 0 &0  & \ldots\\
 -d^{-1}_1 & 0 & \alpha_2+d^{-1}_1+d^{-1}_2 & 0 & -d^{-1}_2 & \ldots\\
  0 & 0 & 0 & 0 & 0  & \ldots\\
  0 & 0 & -d^{-1}_2 & 0 & \alpha_3+d^{-1}_2+d^{-1}_3  & \ldots \\
  \ldots & \ldots &\ldots & \ldots &\ldots   &\ldots \\
\end{array}%
\right).
 \end{equation*}
With respect to the decomposition
$\kH={\widetilde{\kH}}_1\oplus\widetilde{{\kH}}_2$, where
$\widetilde{{\kH}}_1=\spa\{e_{2k-1}\}_{k=1}^\infty$ and
$\widetilde{{\kH}}_2=\spa\{e_{2k}\}_{k=1}^\infty$,
 the operator $T$ admits the representation
$T= S\oplus 0_{\kH_2}
$.
 Hence $\kappa_{-}(T)=\kappa_{-}(S)$, and Theorem  \ref{th8}
 completes the proof.
\end{proof}
Using  equality  $\kappa_{-}(L_{X,\alpha})=\kappa_{-}(S)$ and the
following Gerschgorin theorem, we obtain sufficient condition for
$\kappa_{-}(L_{X,\alpha})\geq m$ as well as for the equality
$\kappa_{-}(L_{X,\alpha})=m$ with arbitrary finite $m$.
\begin{thm}(\cite[Theorem 7.2.1]{Lan})\label{Ger}
All eigenvalues of a matrix
$A=(a_{ij})_{i,j=1}^n\in[\mathbb{C}^n]$ are contained in the union
of Gerschgorin's disks
\[
G_i=\{\lambda\in\mathbb{C}:\,|\lambda-a_{ii}|\leq\sum\limits_{i\neq
j}|a_{ij}|\},\quad i\in\{1,..,n\}.
\]
\end{thm}
\begin{thm}\label{suffcond}
Let matrix $S=(s_{ij})_{i,j=1}^\infty$  be  defined  by
\eqref{eq5}. Suppose that
\begin{equation}\label{suff}
\alpha_{k_i}<-2(d_{k_i-1}^{-1}+d_{k_i}^{-1}) \text{\quad for
\quad} k_i\in K:=\{k_i\}_{i=1}^m.
\end{equation}

Then $\kappa_{-}(L_{X,\alpha})\geq m$.
 If, in  addition,
$\alpha_i>0$ for $i\notin K$, then $\kappa_{-}(L_{X,\alpha})= m.$
\end{thm}
\begin{proof}
Consider  two  cases.

$(a)$ Assume that $k_i=i,\,\, k_i\in K$. Denote by $S_m\in
[\mathbb{C}^m]$ submatrix  in the upper left corner of $S$. In
accordance with the minimax principle (see, for  instance,
\cite{Gla})
\begin{equation}\label{sub}
\kappa_-(S_m)\leq \kappa_-(S)=\kappa_-(L_{X,\alpha}).
\end{equation}
Applying Theorem \ref{Ger} to $S_m$ and using \eqref{suff}, we
obtain $\kappa_-(S_m)=m$. Therefore $\kappa_-(L_{X,\alpha})\geq m$
and the first assertion of the theorem holds.

 Further, setting $\alpha_i=0 \text{ for } i\in\{1,..,m\}$,
we obtain a non-negative self-adjoint operator
$\widetilde{L}_{X,\alpha}$. It is obvious that  $L_{X,\alpha}$ is
an $m$-dimensional perturbation of the operator
$\widetilde{L}_{X,\alpha}$. Thus from  the  minimax principle
follows that $\kappa_-(L_{X,\alpha})\leq m$ and,  consequently,
the second assertion of the theorem is satisfied.

 $(b)$ Let $K$ be an arbitrary set consisting of $m$ natural numbers.
 General  case is easily reduced  to the previous one.
  Namely, there exists  unitary transformation $U$ such that
\begin{equation}
\widetilde{S}=U^*SU,\quad U:\,s_{k_ik_i}\rightarrow
\widetilde{s}_{ii},\quad \sum\limits_{j\neq
k_j}|s_{k_ij}|=\sum\limits_{j\neq i}|\widetilde{s}_{ij}|,\,
\,k_i\in K.
\end{equation}
Applying previous reasoning to the matrix $\widetilde{S}$, we
obtain the proof in the general  case.
\end{proof}
\begin{rem}
Arguing as above, it is not difficult to show that
$\kappa_-(L_{X,\alpha})=\infty$ in the case of infinite $m$.
\end{rem}

 Theorem \ref{th12} enables us to obtain an algorithm for
 determination of
 $\kappa_-(L_{X,\alpha})$. Namely, define the sequence
$\gamma=\{\gamma_k\}_{k=1}^\infty$ by 
\begin{gather}\label{eq8''}
\gamma_1:=\alpha_1+d^{-1}_1,\\
(i)\ \text{if}\quad \gamma_k\neq 0,\quad \text{then}\quad
\gamma_{k+1}:=\alpha_{k+1}+d^{-1}_{k+1}+d^{-1}_k-d_k^{-2}\gamma_k^{-1},\quad
k\geq 1;\label{eq8}\\
(ii)\ \text{if}\quad \gamma_k=0,\quad \text{then}\quad
\begin{array}{c}\gamma_{k+1}:=\infty\\
\gamma_{k+2}:=\alpha_{k+2}+d^{-1}_{k+1}+d^{-1}_{k+2}
\end{array},\quad
k\geq 1.\qquad\quad\label{eq8'}
\end{gather}

\begin{thm}\label{th2}
Let the set $X=\{x_k\}_{k=1}^\infty$ 
satisfy \eqref{d_*^*}. Let the operator $L_{X,\alpha}$ be defined
by \eqref{delta_op} and let the sequence
$\gamma=\{\gamma_k\}^\infty_{k=1}$ be defined by
\eqref{eq8''}--\eqref{eq8'}. Then
\begin{equation*}
\kappa_-(L_{X,\alpha})=\kappa_-(\gamma)+N_\infty(\gamma),
\end{equation*}
where $\kappa_-(\gamma)$ and $N_\infty(\gamma)$ are the number of
negative  and  infinite elements, respectively, in the
 sequence $\gamma$.
\end{thm}
\begin{proof}
Consider two cases.

$(a)$ Let $\gamma_1=\alpha_1+d^{-1}_1\neq0$. 
Setting $T_{11}: = \gamma_1I_{\mathbb{C}}$ and applying Proposition \ref{th2'} to the matrix 
\eqref{eq5}, we get
$\kappa_-(S)=\kappa_-(\gamma_1)+\kappa_-(S_2)$, where
\[
S_2:=\left(%
\begin{array}{cccc}
   \gamma_2&  -d^{-1}_2 & 0&\ldots  \\
   -d^{-1}_2 & \alpha_3+d^{-1}_2+d^{-1}_3 & -d^{-1}_3  &\ldots \\
  0 &  -d^{-1}_3& \alpha_4+d^{-1}_3+d^{-1}_4  & \ldots \\
 \ldots  &\ldots&\ldots &  \ldots   \\
\end{array}%
\right).
\]
Further, if $\gamma_2\neq 0$, then we set $T_{11}=\gamma_2
I_{\mathbb{C}}$ and apply Proposition \ref{th2'} to the matrix
$S_2$. Thus if $\gamma_k\neq 0$ for all $k\in\mathbb{N}$, i.e.,
$N_\infty(\gamma)=0$, then we obtain
$\kappa_-(S)=\kappa_-(\gamma)$.

 $(b)$ Assume that $\gamma_1=\alpha_1+d^{-1}_1=0$. Then
$\gamma_2=\infty$ and $\gamma_3=\alpha_3+d^{-1}_2+d^{-1}_3$.
Let 
$
T_{11}:=\left(\begin{array}{cc}
      0 & -d^{-1}_1 \\
      -d^{-1}_1 & \alpha_{2}+d^{-1}_1+d^{-1}_{2}
\end{array}
\right)\in [\mathbb{C}^2] $.  Since  $\det T_{11}=-d^{-2}_1\neq
0$, by Proposition \ref{th2'}, we get $
\kappa_{-}(S)=\kappa_{-}(T_{11})+\kappa_-(S_3)$, where
\begin{equation*}
S_3:=\left(%
\begin{array}{cccc}
   \gamma_3&  -d^{-1}_3 & 0&\ldots  \\
   -d^{-1}_3 & \alpha_4+d^{-1}_3+d^{-1}_4 & -d^{-1}_4  &\ldots \\
  0 &  -d^{-1}_4& \alpha_5+d^{-1}_4+d^{-1}_5  & \ldots \\
 \ldots  &\ldots&\ldots &  \ldots   \\
\end{array}%
\right).
\end{equation*}
Since $\kappa_-(T_{11})=1$, we get
$\kappa_-(S)=N_\infty(\{\gamma_1,\gamma_2\})+\kappa_-(S_3)$.

Proceeding  as above, we  obtain the desired result.
\end{proof}
Following \cite{AlbNiz03a}, consider  continued  fraction
$A_k:=[\alpha_k; d_{k-1},\alpha_{k-1},..,\alpha_1].$ It is easy to
verify  by induction that if $\gamma_k\neq0$ for all
$\gamma_k\in\gamma$, then
\begin{equation}\label{con frac}
\gamma_k=d_k^{-1}+A_k,\quad k\geq 1.
\end{equation}
 Theorem \ref{th2} and  equality \eqref{con frac} yield  the following  result.
\begin{cor}\label{remunus}
The  operator $L_{X,\alpha}$ is  non-negative  if and only if
\begin{equation*}
A_k>-d_k^{-1},\quad k\geq1.
\end{equation*}
\end{cor}
\subsection{ The case of finite number of $\delta$ -- type  interactions}
Setting $\alpha_k=0, \ k>n$,  in
\eqref{delta_op}, we obtain the operator with
$\delta$-interactions on a finite set.
 Using Theorems \ref{th12} and \ref{th2}, we obtain the following description of the negative squares
$\kappa_{-}(L_{X,\alpha})$.
 Namely, define the  sequence
\begin{equation}\label{eq36''}
\widetilde{\gamma}_1:=\alpha_1+d^{-1}_1,
\end{equation}
\begin{equation}
(i) \  \text{if} \ \widetilde{\gamma}_k\neq0, \ \text{then}\ \
 \widetilde{\gamma}_{k+1}:=\left\{
\begin{array}{cc}
   \alpha_{k+1}+d^{-1}_{k+1}+d^{-1}_k-d_k^{-2}\widetilde{\gamma}_k^{-1},&  k\leq n-1,\\
   d^{-1}_{k+1}+d^{-1}_k-d_k^{-2}\widetilde{\gamma}_k^{-1}, &
k\geq n;
  \end{array} \right.\label{eq36}
  \end{equation}
  \begin{equation}\label{eq36'}
(ii) \  \text{if} \ \widetilde{\gamma}_k=0, \ \text{then}\ \
\widetilde{\gamma}_{k+1}:=\infty,\quad
 \widetilde{\gamma}_{k+2}:=\left\{
\begin{array}{cc}
   \alpha_{k+2}+d^{-1}_{k+1}+d^{-1}_{k+2},&  k\leq n-2,\\
   d^{-1}_{k+1}+d^{-1}_{k+2}, &
k\geq n-1.
  \end{array} \right.
\end{equation}
\begin{cor}\label{th2''}
Let $X=\{x_k\}_{k=1}^{n}\subset\mathbb{R}$ be a  finite set. Let also the operator
 $L_{X,\alpha}$ be defined  by \eqref{delta_op} and
let $\widetilde{\gamma}=\{\widetilde{\gamma}_k\}^\infty_{k=1}$ be the
sequence defined by \eqref{eq36''}--\eqref{eq36'}. Then
\begin{equation*}
\kappa_-(L_{X,\alpha})=\kappa_-(\widetilde{\gamma})+N_\infty(\widetilde{\gamma}).
\end{equation*}
%
%
\end{cor}
 Corollary \ref{th2''} has one essential drawback.
To obtain $\kappa_{-}(L_{X,\alpha})$, we must find infinite number of elements $\widetilde{\gamma}_n,\ n\in \mathbb{N}$. But
 it is possible to overcome this 
by treating $L_{X,\alpha}$ as an extension of the minimal operator
with finite deficiency indices. Namely,
define the
matrix $S\in \dC^{n\times n}$,
\begin{equation}\label{eq25}
S=
\left(%
\begin{array}{ccccc}
  \alpha_1+d^{-1}_1 & -d^{-1}_1 & 0 & \dots & 0 \\
 -d^{-1}_1 & \alpha_2+d^{-1}_1+d^{-1}_2 & -d^{-1}_2 & \dots & 0 \\
  0 & -d^{-1}_2 & \alpha_3+d^{-1}_2+d^{-1}_3 & \dots & 0 \\
 \ldots & \ldots & \ldots & \dots & \ldots \\
  0 & 0 & 0 & \dots & \alpha_n+d^{-1}_{n-1} \\
\end{array}%
\right).
\end{equation}
\begin{thm}\label{thS}
Let $X=\{x_k\}_{k=1}^{n}\subset\mathbb{R}$ be a  finite set. Let the operator
 $L_{X,\alpha}$ be defined  by \eqref{delta_op}  and let
$S$ be the matrix \eqref{eq25}. Then
$\kappa_-(L_{X,\alpha})=\kappa_-(S)$.
\end{thm}
\begin{proof}
Consider the operator $L_{\min}$ of the form \eqref{eq26} with
$X=\{x_k\}_{k=1}^n$. Note that $n_{\pm}(L_{\min})=2n$. The
boundary triplet for  $L_{\min}^*$ might be defined by (cf.
\cite[Section III,\S 1]{Gor84})
\begin{gather}
\mathcal{H}=\oplus^{n}_{k=0}\mathcal{H}_k,\quad
\Gamma_0=\oplus^{n}_{k=0}\Gamma_0^k,\quad
\Gamma_1=\oplus^{n}_{k=0}\Gamma_1^k,\quad \text{where} \label{btII_1}\\
\mathcal{H}_0=\mathbb{C},\quad\Gamma_0^0f=-f(x_1-),\quad
\Gamma_1^0f=f'(x_1-),\\
\mathcal{H}_k=\mathbb{C}^2,\,\Gamma_0^kf=\binom{f(x_k+)}{-f(x_{k+1}-)},
\,\Gamma_1^kf=\binom{f'(x_k+)} {f'(x_{k+1}-)}, k\in\{1,..,n-1\},\\
\mathcal{H}_n=\mathbb{C},\quad \Gamma_0^nf=f(x_n+),\quad
\Gamma_1^nf=f'(x_n+).\label{btII_4}
 \end{gather}
  The corresponding Weyl function $M(\lambda)$ is
\begin{equation*}
M(\lambda)=\oplus_{k=0}^nM_k(\lambda),\quad
M_0(\lambda)=M_n(\lambda)=i\sqrt{\lambda} \text{\quad and }
\end{equation*}
$M_k(\lambda)$ for  $k\in\{1,..,n-1\}$ is given  by \eqref{weyl}.

Using \eqref{btII_1}--\eqref{btII_4}, we obtain a description of
$\dom(L_{X,\alpha})$ in the  form \eqref{eq80},
where
\begin{equation*}
C=
 \left(%
   \begin{array}{ccccc}
   0 & \alpha_1 & \ldots &0&0  \\
   1 & 1 &  \ldots  &0&0  \\
\ldots &\ldots & \ldots &\ldots &\ldots \\
0&0&\ldots &0&\alpha_n\\
0&0&\ldots &1&1\\
  \end{array}%
\right), \,
D=\left(%
   \begin{array}{ccccc}
  -1 & 1 &  \ldots & 0&  0 \\
   0&0& \ldots & 0 &  0    \\
 \ldots & \ldots&\ldots& \ldots & \ldots\\
 0&0&\ldots&-1&1\\
 0&0&\ldots&0&0\\
  \end{array}
\right).
\end{equation*}
Further, it easy to  verify  that the matrix $T=CD^*-DM(0)D^*$ has
the form
\begin{equation*}
T=
\left(%
\begin{array}{cccccc}
  \alpha_1+d^{-1}_1 & 0 & -d^{-1}_1 & \ldots & 0&0  \\
  0 & 0 & 0 & \ldots &0  & 0\\
 -d^{-1}_1 & 0 & \alpha_2+d^{-1}_1+d^{-1}_2 & \ldots & 0 & 0\\
  \ldots& \ldots & \ldots & \ldots & \ldots  & \ldots\\
  0 & 0 & 0 & \ldots & 0  & 0 \\
  0 & 0 &0 & \ldots &0   & \alpha_n+d^{-1}_n\\
\end{array}%
\right).
 \end{equation*}
 Arguing as in the proof of Theorem \ref{th12}, we complete the proof.
\end{proof}
 Define the sequence $\gamma=\{\gamma_k\}^n_{k=1}$ as
follows
\begin{gather}\label{eq58''}
\gamma_1:=\alpha_1+d^{-1}_1,\\
(i)\ \text{if}\ \gamma_k\neq 0,\ \text{then} \quad
\gamma_{k+1}:=\left\{
\begin{array}{cc}
   \alpha_{k+1}+d^{-1}_{k+1}+d^{-1}_k-d_k^{-2}\gamma_k^{-1},&  k\leq n-2\\
   \alpha_k+d^{-1}_{k-1}-d_{k-1}^{-2}\gamma_{k-1}^{-1}, &
k= n-1
  \end{array} \right. ;\label{eq58}\\
(ii)\ \text{if} \ \gamma_k=0,\ \ \text{then} \
\begin{array}{c}
\gamma_{k+1}:=\infty,\quad k\in\{1,..,n-1\}, \quad \\
\gamma_{k+2}:=\alpha_{k+2}+d^{-1}_{k+1}+d^{-1}_{k+2},\quad
k\in\{1,..,n-2\}.
\end{array}
\label{eq58'}
 \end{gather}

\begin{thm}\label{th13}
Assume $X=\{x_k\}_{k=1}^{n}$. 
Let the operator $L_{X,\alpha}$ be defined  by \eqref{delta_op} and
 let the  sequence $\gamma=\{\gamma_k\}^n_{k=1}$ be  defined by
\eqref{eq58''}--\eqref{eq58'}. Then
\[
\kappa_-(L_{X,\alpha})=\kappa_-(\gamma)+N_\infty(\gamma).
\]
\end{thm}
We omit the proof since it is analogous to that of Theorem \ref{th2}.
\begin{prop}
Corollary \ref{th2''} and Theorem \ref{th13} are equivalent, i.e.,
\[
\kappa_-(\widetilde{\gamma})+N_\infty(\widetilde{\gamma})=\kappa_-(\gamma)+N_\infty(\gamma),
\]
where $\widetilde{\gamma}=\{\widetilde{\gamma}_n\}_{n=1}^\infty$ and $\gamma=\{\gamma_n\}_{n=1}^\infty$ are defined by \eqref{eq36''}--\eqref{eq36'} and \eqref{eq58''}--\eqref{eq58'}, respectively.
\end{prop}
\begin{proof}
Since $\widetilde{\gamma_k}=\gamma_k$ for $k<n$, it suffices to verify that
\begin{equation}\label{3.20}
\kappa_{-}(\gamma_n)+N_\infty(\gamma_n)=\kappa_{-}(\{\widetilde{\gamma_k}\}_{k=n}^\infty)+N_\infty(\{\widetilde{\gamma_k}\}_{k=n}^\infty).
\end{equation}



First, assume that $\widetilde{\gamma}_m<0$ for some $m\geq n$. Then, by \eqref{eq36}, $\widetilde{\gamma}_{m+1}>d_{m+1}^{-1}$ and hence
\begin{equation}\label{ind}
\widetilde{\gamma}_{k}\geq d^{-1}_{k}, \text{ for all } k>m.
\end{equation}
 The latter also yields that in this case $\kappa_{-}(\{\widetilde{\gamma_k}\}_{k=n}^\infty)\leq 1$ and $N_\infty(\{\widetilde{\gamma_k}\}_{k=m}^\infty)=0$.

Further, if $\widetilde{\gamma}_m=0$ for some $m> n$, then, by \eqref{eq36'},
$\widetilde{\gamma}_{m+1}=\infty$ and
$\widetilde{\gamma}_{m+i}=d^{-1}_{m+i}+d^{-1}_{m+i-1}$, $i\geq 2$. Therefore, $N_\infty(\{\widetilde{\gamma_k}\}_{k=n}^\infty)\leq 1$ and $\kappa_{-}(\{\widetilde{\gamma_k}\}_{k=m}^\infty)=0$.


Consider three cases.

 $\textbf{(a)}$ Let $\gamma_n\geq0$. Combining \eqref{eq58} with \eqref{eq36}, we get  $\widetilde{\gamma}_n=\gamma_n+d_n^{-1}\geq d_n^{-1}$. By \eqref{eq36}, $\widetilde{\gamma}_k$ satisfies  \eqref{ind} with $m=n$,
and hence \eqref{3.20} clearly holds.

$\textbf{(b)}$ Let $\gamma_n=\infty$. Then
$\gamma_{n-1}=\widetilde{\gamma}_{n-1}=0$ and
$\widetilde{\gamma}_n=\infty$. Thus $\widetilde{\gamma}_k$ satisfies  \eqref{ind} for all $k\geq n$, and hence \eqref{3.20} holds.

$\textbf{(c)}$ Assume now that $\gamma_n<0$.

If $\widetilde{\gamma}_k = 0$ for some $k\geq n$, then arguing as above we arrive at \eqref{3.20}.

Suppose that $\widetilde{\gamma}_k\neq0$, $k\geq n$. To prove
\eqref{3.20} it suffices to show that $\widetilde{\gamma}_{n+i}<0$
for some $i>0$. Assume the converse, i.e.,
$\widetilde{\gamma}_{k}>0$ for all $k\geq n$. Denote
$\xi_k:=d^{-1}_k-\widetilde{\gamma}_k$, $k\geq n$. Clearly,
$\xi_k<d_k^{-1}$ for all $k\geq n$. Note  that,
by \eqref{eq36},  $\xi_{k+1}=d^{-1}_k\xi_k(d^{-1}_k-\xi_k)^{-1}$.
Further, the  inequality  $0<\xi_n<d^{-1}_n$ holds since
$\gamma_n<0$. Moreover, $0<d^{-1}_n-\xi_n<d^{-1}_n<d_*^{-1}$ (see
\eqref{d_*^*}) and hence
\begin{equation*}
\xi_{n+1}=d^{-1}_n\xi_n(d^{-1}_n-\xi_n)^{-1}=\xi_n+\xi_n^2(d^{-1}_n-\xi_n)^{-1}>\xi_n+\xi_n^2d_*.
\end{equation*}
Similarly, 
$\xi_{n+1}<d^{-1}_{n+1}\leq d_*^{-1}$ yields 
\[
\xi_{n+2}>\xi_{n+1}+\xi_{n+1}^2d_*>\xi_n+\xi_n^2d_*+\xi_n^2d_*=\xi_n+2\
\xi_n^2d_*.
\]
Therefore, we get $\xi_{n+i}>\xi_n+i\ \xi_n^2d_*$, $i\in
\mathbb{N}$. Hence there exists $i_0\in \mathbb{N}$ such that
\[
\xi_n+i \xi_n^2d_*>d_*^{-1}>d_{n+i}^{-1},\qquad i\geq i_0.
\]
Therefore we get $\xi_{n+i_0}>d^{-1}_{n+i_0}$ and consequently
$\widetilde{\gamma}_{n+i_0}<0$. This contradiction comletes the
proof of \eqref{3.20}.

Combining $\textbf{(a)}$, $\textbf{(b)}$, and $\textbf{(c)}$, we
arrive at the desired result.
\end{proof}

\begin{rem}\label{AN appr}
In \cite{AlbNiz03}, S. Albeverio and L. P. Nizhnik obtained
another description of $\kappa_{-}(L_{X,\alpha})$. Namely, define the function $\varphi$ as a solution
of the problem
\begin{gather}\label{eq6} \varphi''(x)=0, \quad
x\notin X,\qquad \varphi(x)\equiv 1,\quad  x<x_1,\quad
\text{and}\\
\varphi(x_k+)=\varphi(x_k-), \quad
                     \varphi'(x_k+)-\varphi'(x_k-)=\alpha_k\varphi(x_k) \quad\text{for}\quad x_k\in X.
\label{eq6'}\end{gather}
Theorem 3  from \cite{AlbNiz03} states that \emph{
$\kappa_{-}(L_{X,\alpha})$ 
equals the  signature of the sequence
\begin{equation}
(\varphi(x_1),\,\varphi(x_2),...,\,\varphi(x_n),\,(1+\alpha_n
d_{n-1})\varphi(x_n)-\varphi(x_{n-1})).
\end{equation}}
 Note that this result  may be deduced from Theorem
\ref{thS} and vise versa. Namely, let $\Delta_k$ be a $k$-th order
leading principle minor of the matrix $S$ defined by \eqref{eq25}.
Then one can check that
\begin{equation*}\label{eq15}
\Delta_k=\frac{\varphi(x_{k+1})}{d_{k-1}\cdot...\cdot d_1},\qquad
k\in\{1,..,n\}.
\end{equation*}
 \end{rem}

 \section{Operators with $\delta'$-interactions}


The main result of this Section is the following theorem.
\begin{thm}\label{th30}
Let $X=\{x_k\}_{k=-\infty}^\infty$ be a discrete subset of
$\mathbb{R}$ satisfying \eqref{d_*^*},  $L_{X,\beta}$  the
operator defined by \eqref{delta'_op}, and
$\beta=\{\beta_k\}_{k=-\infty}^\infty\subset\mathbb{R}$. Then
 $\kappa_{-}(L_{X,\beta})=\kappa_{-}(\beta)$.
   \end{thm}

  \begin{proof}
   We divide the proof into several steps.

   \textbf{(a)}
Consider the minimal operator  $L_{\min}$ \eqref{eq26}. Since $X$
satisfies \eqref{d_*^*}, we can choose the boundary triplet
$\Pi=\{\kH,\gG_0,\gG_1\}$  for $L_{\min}^*$ as follows \cite[Lemma
1]{Koch}
\begin{equation}\label{eq27b}
\mathcal{H}=\oplus^{\infty}_{k=-\infty}\mathcal{H}_k,\quad
\Gamma_0f=\oplus^{\infty}_{k=-\infty}\Gamma_0^kf,\quad
\Gamma_1f=\oplus^{\infty}_{k=-\infty}\Gamma_1^kf,
\end{equation}
 where $\Pi_k=\{\mathcal{H}_k,\Gamma_0^k,\Gamma_1^k\}$, $k\in\mathbb{Z}$, is given by \eqref{eq27''}.

  The corresponding Weyl function is $M(\lambda)=\oplus^{\infty}_{k=-\infty} M_{k}(\lambda)$ with $M_k(\lambda)$ defined
  by \eqref{weyl}.  

The domain  of the operator $L_{X,\beta}$ admits the
representation
$\dom(L_{X,\beta})=\dom(L_{\min}^*)\upharpoonright\ker(D\Gamma_1-C\Gamma_0)$
with $D$ and $C$ determined, respectively, by
\begin{equation*}
\left(%
\begin{array}{ccccccc}
 \ldots&  \ldots& \ldots & \ldots&\ldots  &\ldots\\
  \ldots &  -1& 1 & 0 & 0 & \ldots \\
  \ldots &0  & \ulcorner \beta_0 &\ 0\ \urcorner &0& \ldots \\
    \ldots& 0 &\llcorner\ 0\  & -1 \lrcorner&1 &\ldots \\
 \ldots &   0 & 0 & 0 & \beta_1&\ldots \\
 \ldots & \ldots& \ldots & \ldots&\ldots &\ldots\\
 \end{array}%
\right)\text{and}
\left(%
\begin{array}{ccccccc}
 \ldots&  \ldots& \ldots & \ldots& \ldots &\ldots\\
  \ldots &  0& 0 & 0 & 0 & \ldots \\
  \ldots &1  & \ulcorner 1 & 0 \urcorner &0& \ldots \\
    \ldots& 0 &\llcorner 0  & 0 \lrcorner&0 &\ldots \\
 \ldots &   0 & 0 & 1 & 1&\ldots \\
 \ldots & \ldots&\ldots  &\ldots &\ldots &\ldots\\
 \end{array}%
\right).
\end{equation*}

  Arguing as in the  proof of Theorem \ref{th12}, we  assume that $D$
is bounded. After straightforward calculations we get the operator
$T=CD^*-DM(0)D^*$,
\begin{equation}\label{4.4}
 \left(%
\begin{array}{ccccccccc}
 \ldots&  \ldots & \ldots &  \ldots& \ldots &\ldots& \ldots \\
 \ldots  & \beta_{0}d^{-1}_{-1} & -d^{-1}_{-1} & 0 & 0 &0 & \ldots\\
 \ldots& \ulcorner\beta_{0}+{\beta_{0}}^2d^{-1}_{-1} & -\beta_{0}d^{-1}_{-1}\urcorner & 0 & 0 &0 & \ldots\\
  \ldots &\llcorner -\beta_{0}d^{-1}_{-1} & d^{-1}_{-1}+d^{-1}_{0}\lrcorner & \beta_1d^{-1}_{0} & -d^{-1}_{0} & 0&\ldots \\
   \ldots  & 0 &\beta_1d^{-1}_{0} & \ulcorner\beta_1+{\beta_1}^2d^{-1}_{0} & -\beta_1d^{-1}_{0}\urcorner & 0&\ldots \\
   \ldots  & 0 &  -d^{-1}_{0} & \llcorner-\beta_1d^{-1}_{0}  & d^{-1}_{0}+d^{-1}_1\lrcorner & \beta_2d^{-1}_1&\ldots \\
\ldots&\ldots   &   \ldots & \ldots & \ldots &\ldots& \ldots   \\
\end{array}%
\right).
\end{equation}
By Theorem \ref{th8}, $\kappa_{-}(L_{X,\beta})=\kappa_{-}(T)$.

\textbf{(b)} Note that the matrix \eqref{4.4}  admits the
 representation $T=A+B$, where
\[
A=\sum_{k=-\infty}^\infty\beta_k(\cdot,e_{2k-1})e_{2k-1},\quad B=
\sum\limits^\infty_{k=-\infty}d_{k-1}^{-1}(\cdot,\mathrm{\mathbf{b}}_k)\mathrm{\mathbf{b}}_k, 
\]
with $\mathrm{\mathbf{b}}_k:=e_{2k-1}+\beta_{k}e_{2k}-e_{2k+1}$.
Since  $d^{-1}_k>0$, one gets
\begin{equation}\label{eq39}
\kappa_{-}(T)\ \leq \ \kappa_{-}(A). 
\end{equation}

 \textbf{(c)}
  Let $s\in \dZ_{-}\cup\{0\}$ and $r\in \dN$.
Consider the matrix
$T_{s,r}\in
 [\mathbb{C}^{2(r-s)+1}]$,
\begin{gather*}
T_{s,r}:=A_{s,r}+B_{s,r}, \quad\text{where}\quad A_{s,r}=\sum_{k=s-1}^{r+1}\beta_k(\cdot,e_{2k})e_{2k} \text{ and}\\
%
B_{s,r}=\sum_{k=s+1}^{r-1}d_{k-1}^{-1}(\cdot,\mathrm{\mathbf{b}}_k)\mathrm{\mathbf{b}}_k+
d^{-1}_{s-1}(\cdot,\mathrm{\mathbf{y}}_s)\mathrm{\mathbf{y}}_s
+d^{-1}_{r-1}(\cdot,\mathrm{\mathbf{x}}_r)\mathrm{\mathbf{x}}_r,
\end{gather*}
with $\mathrm{\mathbf{y}}_s:=\beta_se_{2s-2}+e_{2s-1}$ and
$\mathrm{\mathbf{x}}_r:=e_{2r+1}+\beta_re_{2r+2}$.

 It is clear that $
\ran(A_{s,r})\cap\ran\left(B_{s,r}-d^{-1}_{s-1}(\cdot,\mathrm{\mathbf{y}}_s)\mathrm{\mathbf{y}}_s\right)=\{0\}
$ and hence
\[
\kappa_{-}(T_{s,r}-d^{-1}_{s-1}(\cdot,\mathrm{\mathbf{y}}_s)\mathrm{\mathbf{y}}_s)
=\kappa_{-}(A_{s,r}).
\]
According to the choice  of the matrix $T_{s,r}$, we get
\begin{equation}\label{eq40}
\kappa_{-}(T)\geq
\kappa_{-}(T_{s,r}-d^{-1}_{s-1}(\cdot,\mathrm{\mathbf{y}}_s)\mathrm{\mathbf{y}}_s)
-\rank(d^{-1}_{s-1}(\cdot,\mathrm{\mathbf{y}}_s)\mathrm{\mathbf{y}}_s)=\kappa_{-}(A_{s,r})-1.
\end{equation}
Combining \eqref{eq39} with \eqref{eq40}, we obtain
$  \kappa_{-}(A)-1\leq\kappa_{-}(T)\leq \kappa_{-}(A)$.

If $\kappa_-(\beta)=\infty$, then \eqref{eq39} yields $\kappa_{-}(T)=\infty$. Therefore, $\kappa_{-}(L_{X,\beta})=\infty$ and theorem is proven in the case $\kappa_-(\beta)=\infty$.

 \textbf{(d)} Assume now that  $\kappa_{-}(\beta)=m<\infty$.
 Let us show that
 \begin{equation}\label{eq70}
\det(T_{s,r})=\bigl(x_{r}-x_{s-1} +\sum\limits^{r}_{k=s}
\beta_k\bigr)\prod\limits^{r}_{k=s}d^{-1}_{k-1}\beta_k.
 \end{equation}
 For $s=0,\,r=1$ equality \eqref{eq70} is obvious. Suppose that \eqref{eq70} holds with  $s+1<0$ and
 $r>1$.
 Note that the  second  row $\mathrm{\mathbf{t}}_2$  of the  matrix $T_{s,r}$ admits  a decomposition
 \begin{gather*}
 \mathrm{\mathbf{t}}_2=\mathrm{\mathbf{t}}_2^1+\mathrm{\mathbf{t}}_2^2,
 \text{ where}\\
 \mathrm{\mathbf{t}}_2^1=\begin{pmatrix}
  -\beta_{s}d^{-1}_{s-1} & d^{-1}_{s-1} & 0 & \ldots & 0\\
\end{pmatrix} \text{ and } \mathrm{\mathbf{t}}_2^2=\begin{pmatrix}
 0 & d^{-1}_{s} &  \beta_{s+1}d^{-1}_{s} & \ldots & 0\\
\end{pmatrix}.
\end{gather*}
 Then $\det(T_{s,r})=\det(T^1_{s,r})+\det(T^2_{s,r})$, where
 $T^1_{s,r}$ and $T^2_{s,r}$
   are  matrices  obtained by replacement of  $\mathrm{\mathbf{t}}_2$  in  $T_{s,r}$ by $\mathrm{\mathbf{t}}_2^1$
   and  $\mathrm{\mathbf{t}}_2^2$, respectively.
 Adding  to the  first  row    of
$T^1_{s,r}$   the second row multiplied  by $\beta_s$, we arrive
at the equality
\begin{equation}\label{eq71}
\det(T^1_{s,r})=\beta_sd^{-1}_{s-1}\det(T_{s+1,r}).
\end{equation}
It is easily seen  that $\det(T^2_{s,r})=(\beta_{s}+
\beta_{s}^2d^{-1}_{s-1})\det(T_{2(r-s)}^2)$,
 where $T_{2(r-s)}^2\in
[\mathbb{C}^{2(r-s)}]$ is an algebraic complement of $\beta_{s}+
\beta_{s}^2d^{-1}_{s-1}$. Adding to the second  row of
$T_{2(r-s)}^2$  the  first row multiplied by $-\beta_{s+1}$ and to
the  third  row the  first row, we get
$\det(T^2_{s,r})=(\beta_{s}+
\beta_{s}^2d^{-1}_{s-1})\beta_{s+1}d_{s}^{-1}\det(T_{2(r-s-1)}^2).$
 Proceeding  analogously, one, finally, obtains
\begin{equation}\label{eq72}
 \det(T^2_{s,r})=(\beta_{s}+
\beta_{s}^2d^{-1}_{s-1})\bigl(\prod\limits^{r}_{k=s+1}d^{-1}_{k-1}\beta_k\bigr).
\end{equation}
Combining \eqref{eq71} and \eqref{eq72}, we arrive at
\eqref{eq70}.

 Since $\kappa_-(\beta)<\infty$ and the difference $(x_{r}-x_{s-1})=\sum\limits_{k=s}^{r}d_k$ is
 unbounded as either  $-s$ or $r$ tends to infinity,  for sufficiently large
 $s$ and $r$  we have
$ \bigl(x_{r}-x_{s-1} +\sum\limits^{r}_{k=s} \beta_k\bigr)>0$ and
$\beta_k>0$ for $k<s$ or $k>r$.
 Hence
\begin{equation*}
\sgn(\det(T_{s,r}))=\sgn\bigl(\prod\limits^{r}_{k=s}d^{-1}_{k-1}\beta_k\bigr)=(-1)^{\kappa_-(\beta)}=(-1)^m\neq(-1)^{m-1}.
\end{equation*}
Therefore $\kappa_{-}(T_{s,r})= m$ and, by \eqref{eq39}, we,
finally, get
\begin{equation*}
m=\kappa_{-}(T_{s,r})\leq \kappa_{-}(T)\leq \kappa_{-}(A)=m,
\end{equation*}
i.e., $\kappa_{-}(T)=m$. The proof is completed.
\end{proof}
\subsection*{Acknowledgment}
The authors are grateful to M.M. Malamud for posing of the problem
and  permanent  attention to our  work. We are especially indebted
to A.S. Kostenko for carefully  reading of the preliminary version
of the manuscript and constructive remarks. 
 We would  also  like  to thank the referee for useful
remarks  regarding improving of the exposition.

\end{document}